\documentstyle[12pt,leqno]{article}

\input amssym.def
\input amssym

\def\twelveamsfonts{
 \font\twelvemsa=msam10 scaled 1200
 \font\eightmsa=msam8
 \font\sixmsa=msam6
 \font\twelvemsb=msbm10 scaled 1200
 \font\eightmsb=msbm8
 \font\sixmsb=msbm6
 \font\twelvembi=cmmib10 scaled 1200
 \font\eightmbi=cmmib8
 \font\sixmbi=cmmib6
 \textfont\msafam\twelvemsa
 \scriptfont\msafam\eightmsa
 \scriptscriptfont\msafam\sixmsa
 \textfont\msbfam\twelvemsb
 \scriptfont\msbfam\eightmsb
 \scriptscriptfont\msbfam\sixmsb}

\ifnum\ht0<\ht2\let\twelveamsfonts\relax\fi
\twelveamsfonts

\let\bls\baselineskip \let\nt\noindent
\let\vp\vphantom 
  
\def\vsk#1>{\vskip#1\bls} \def\vv#1>{\vadjust{\vsk#1>}}
\def\,{\relax\ifmmode\mskip\thinmuskip\relax\else\kern.16667em\fi}
\def\;{\relax\ifmmode\mskip\thickmuskip\relax\else\kern.27777em\fi}
\def\!{\relax\ifmmode\mskip-\thinmuskip\relax\else\kern-.16667em\fi}
\def\&{.\kern.1em} \def\itl#1{{\it #1\/}} 
\def\ftext#1{{\let\thefootnote\relax\footnotetext{\vsk-.8>\noindent #1}}}

\let\geq\geqslant

\let\leq\leqslant

\def\bea{\begin{eqnarray}}
\def\ena{\end{eqnarray}}
\def\beq{\begin{equation}}
\def\eeq{\end{equation}}
\def\no{\nonumber}
\def\qed{\quad$\square$}

\def\proof{\noindent{\it Proof.}\quad}
\def\rem{\noindent{\it Remark.}\quad}
\newtheorem{thm}{Theorem}[section]
\newtheorem{prop}[thm]{Proposition}

\makeatletter
\newbox\p@b@ld
\def\poorbold#1{\setbox\p@b@ld\hbox{#1}\kern-.01em\copy\p@b@ld\kern-\wd\p@b@ld
 \kern.02em\copy\p@b@ld\kern-\wd\p@b@ld\kern-.012em\raise.02em\box\p@b@ld}

\@addtoreset{thm}{section}
\@addtoreset{equation}{section}
\makeatother

\begin{document}

\begin{center}
\vp1
{\Large \bf The $q$-twisted Cohomology and
\vsk.3>
The $q$-hypergeometric  Function at $|q|=1$
}
\vsk2>
{Yoshihiro Takeyama$^{\,\diamond}$}
\ftext{$^{\diamond\,}$Research Fellow of the Japan Society for
the Promotion of Science.
\\ E-mail: ninihuni@kurims.kyoto-u.ac.jp
}
\vsk1.5>
{\it Research Institute for Mathematical Sciences,
Kyoto University, Kyoto 6068502, Japan}
\end{center}
\vsk1.75>

{\narrower\nt
{\bf Abstract.}\enspace
We construct the $q$-twisted cohomology associated with the $q$-multiplicative function 
of Jordan-Pochhammer type at $|q|=1$. 
In this framework, we prove the Heine's relations and a connection formula 
for the $q$-hypergeometric function of the Barnes type. 
We also prove an orthogonality relation of
the $q$-little Jacobi polynomials at $|q|=1$.
\vsk1.4>}
\vsk0>
\thispagestyle{empty}

\section{Introduction}
In this paper we construct the $q$-twisted cohomology at $|q|=1$ in Jordan-Pochhammer case 
and prove some properties of the $q$-hypergeometric function at $|q|=1$ defined in \cite{NU}.

The basic hypergeometric function with $0<|q|<1$ \cite{GR} is represented in terms of a Jackson integral. 
In \cite{A1, AK}, a formulation of Jackson integrals is given. 
Namely, for a $q$-multiplicative function defined by means of $q$-version of Sato's $b$-functions \cite{SSM}, 
the $q$-twisted cohomology is defined. 
In this approach, Jackson integrals can be regarded as a pairing between this cohomology and $q$-cycles. 

We consider the case that $|q|=1$ and $q$ is not a root of unity. 
Then the structure of $b$-functions is the same as in the case of $0<|q|<1$ 
and associated $q$-multiplicative function can be constructed in terms of the double sine function \cite{B}. 
The problem is to define a suitable integral which is a certain generalization of Jackson integrals 
to the case of $|q|=1$. 

In \cite{MT}, a family of solutions to the quantum Knizhnik-Zamolodchikov ($q$KZ) equation at $|q|=1$ 
was constructed. 
The solution is represented in terms of a pairing between two functional spaces, 
which is called the hypergeometric pairing. 
The hypergeometric pairing was defined by Tarasov and Varchenko 
in the study of the rational $q$KZ \cite{TV1} and the trigonometric $q$KZ for $0<|q|<1$ \cite{TV2}. 
In the trigonometric case, the hypergeometric pairing is a pairing 
between a space of trigonometric functions and that of elliptic functions. 
It is given by an integral over a closed contour 
with a kernel function defined by a product of the infinite product with step $q$. 
This kernel function is the $q$-multiplicative funtion mentioned above. 
By taking residues, we can represent this integral in terms of a Jackson integral.  

In the case of $|q|=1$, the hypergeometric pairing is a pairing between two spaces of trigonometric functions 
which depend on two respective values of a deformation parameter
\bea
q=e^{2 \pi i \omega} \quad {\rm and} \quad Q=e^{\frac{2 \pi i}{\omega}}.
\no
\ena 
Hence the pairing induces certain duality of two functional spaces at $|q|=1$. 
This type of duality has appeared in mathmatical physics: for example, 
modular double of quantum group \cite{F} and matrix elements in quantum Toda chain \cite{S}.

In this paper we define an integral associated with the $q$-multiplicative function 
of Jordan-Pochhammer type at $|q|=1$ in a similar way to \cite{MT}. 
This integral is regarded as a pairing between two functional spaces which depend on $q$ and $Q$, respectively. 
Then we can define a cohomology of these spaces associated with this integral. 
In this way we get the $q$-twisted cohomology at $|q|=1$. 
{}From this point of view, we can prove some relations satisfied by the $q$-hypergeometric function 
of the Barnes type at $|q|=1$.

The plan of the paper is as follows. 
In Section 2, we recall the result of $q$-analogue of $b$-funcions following \cite{A2} and
define a $q$-multiplicative function at $|q|=1$. 
In Section 3, we construct the $q$-twisted cohomology associated with the $q$-multiplicative function 
of Jordan-Pochhammer type. 
In this case we can find a basis of the cohomology. 
In order to prove linear independence, we write down the formula for a determinant, see (\ref{7.9}). 
This formula gives us the $q$-Beta integral formula at $|q|=1$ in a special case. 
In Section 4, we prove two properties of the $q$-hypergeometric function at $|q|=1$: 
Heine's relations and a connection formula. 
Section 5 is additional one. 
We discuss the $q$-little Jacobi polynomials and their orthogonality 
with respect to the kernel of the $q$-Beta integral at $|q|=1$ given in Section 3. 

\section{The $q$-multiplicative function at $|q|=1$}
Let $q$ be a nonzero complex number. 
In this paper, we consider the case that $|q|=1$ and $q$ is not a root of unity. 
We put $q=e^{2 \pi i \omega} \,  (\omega >0,  \omega \not\in {\Bbb Q}\, ).$

Let $L$ be an $l$ dimensional integer lattice in ${\Bbb C}^{l}$:
\bea
L:=\{ \chi =(\chi_{1}, \cdots , \chi_{l}) | \chi_{j} \in {\Bbb Z}, j=1, \cdots , l \}
\in {\Bbb C}\,^{l}.
\ena
For a set of nonzero rational functions $\{b_{\chi}(t)\}_{\chi \in L}$, where
\bea
b_{\chi}(t)=b_{\chi}(t_{1}, \cdots , t_{l}) \in {\Bbb C}\,(t_{1}, \cdots , t_{l})^{\times}, 
\ena
we consider the following system of difference equations:
\bea
\Phi(z+\chi)=b_{\chi}(t)\Phi(z), \quad (\chi \in L),
\label{1}
\ena
where $z=(z_{1}, \cdots , z_{l}) \in {\Bbb C}^{l}$ and
\bea
t=(t_{1}, \cdots , t_{l}):=(e^{2 \pi i \omega z_{1}}, \cdots , e^{2 \pi i \omega z_{l}}).
\ena
The compatibility condition of (\ref{1}) implies
\bea
& & b_{0}(t)=1, \label{2} \\
& & b_{\chi + \chi'}(t)=b_{\chi}(t)b_{\chi'}(q^{\chi}\!\cdot\! t) \quad {\rm for \, any} \, \chi, \chi' \in L,
\label{3} 
\ena
where $ q^{\chi}\!\cdot\! t=(q^{\chi_{1}}t_{1}, \cdots , q^{\chi_{l}}t_{l})$.
The conditions (\ref{2}) and (\ref{3}) mean that the set $\{b_{\chi}(t)\}_{\chi \in L}$ defines a 
{\it 1-cocycle}. 
A set $\{b_{\chi}(t)\}_{\chi \in L}$ is said to be a {\it 1-coboundary} if and only if 
there exists a nonzero rational funcion $\varphi (t)$ such that
\bea
b_{\chi}(t)=\frac{\varphi(q^{\chi}\cdot t)}{\varphi(t)}  \quad {\rm for \, any} \, \chi \in L.
\ena

Let us consider the quotient
$H^{1}:=$\{1-cocycles\}$/$\{1-coboundaries\}.
$H^{1}$ has a mutiplicative group structure.
For $\mu \in L^{*}:={\rm Hom}_{\Bbb Z}(L, {\Bbb Z})$, we put
\bea
\mu(m)=\mu(0, \cdots , \stackrel{m-{\rm th}}{1}, \cdots , 0) \in {\Bbb Z}
\ena
and
\bea
t^{\mu}=t_{1}^{\mu(1)}\cdots t_{l}^{\mu(l)}.
\ena
Note that $\mu(\chi)=\sum_{m=1}^{l}\mu(m)\chi_{m}$ for $\chi=(\chi_{1}, \cdots , \chi_{l}) \in L$.

Then the following result holds.
\begin{prop}
$H^{1}$ is represented by cocycles of the following form:
\bea
b_{\chi}(t)=a_{\chi}\prod_{\nu=0}^{\mu_{0}(\chi)-1}(q^{\nu}t^{\mu_{0}})
                    \frac{\displaystyle \prod_{j=1}^{k}(q^{\gamma_{j}}t^{\mu_{j}}; q)_{\mu_{j}(\chi)}}
                         {\displaystyle \prod_{j=1}^{k'}(q^{\gamma_{j}'}t^{\mu_{j}'}; q)_{\mu_{j}'(\chi)}}
\label{4}
\ena
for $\mu_{0}, \mu_{j}, \mu_{j}' \in L^{*}$ and $\gamma_{j}, \gamma_{j}' \in {\Bbb C}\, $. 
Here $\{a_{\!\chi}\}_{\chi \in L}$ is a set of nonzero constants satisfying $a_{\chi + \chi'}=a_{\chi}a_{\chi'}$
for any $\chi, \chi' \in L$, and
\bea
(x; q)_{n}:=\left\{ \begin{array}{l} \displaystyle \prod_{j=1}^{n-1}(1-xq^{j}), \quad {\rm for} \, n \geq 0, \\
                                  \displaystyle \prod_{j=1}^{-n}(1-xq^{-j})^{-1}, \quad {\rm for} \, n<0.
                 \end{array} \right.
\ena
The expression (\ref{4}) is not unique.
\end{prop}

This result is a $q$-analogue of Sato's result \cite{SSM}
and was stated by Aomoto in \cite{A2}. 
In \cite{A2}, Aomoto stated this result in the case of $0<|q|<1$. 
However, it can be checked that the proposition holds unless $q$ is a root of unity.

Let us find a solution $\Phi(z)$ to (\ref{1}) 
for the 1-cocycle $\{b_{\chi}(t)\}_{\chi \in L}$ given by (\ref{4}).
Following \cite{N}, we set
\bea
\langle x \rangle:=
\exp{(\frac{\pi i}{2} \left( (1+\omega)x-\omega x^{2} \right))}S_{2}(x | 1, \frac{1}{\omega}), 
\ena
where $S_{2}(x)$ is the double sine function.
We refer the reader to \cite{JM} for the double sine function. 
Moreover, we define a funcion $\sigma(x)$ by
\bea
\sigma(x):=\exp{(\pi i\left( (1+\omega)x-\omega x^{2} \right))}
         =\langle x \rangle \langle 1+\frac{1}{\omega}-x \rangle.
\label{4.1}
\ena
These functions satisfy
\bea
\frac{\langle x+1 \rangle}{\langle x \rangle}=\frac{1}{1-e^{2 \pi i \omega x}}, \quad
\frac{\sigma(x+1)}{\sigma(x)}=-e^{-2\pi i \omega x}.
\label{4.9}
\ena
For $\mu \in L^{*}$, we set
\bea
\mu(z):=\sum_{m=1}^{l}\mu(m)z_{m}.
\ena
Then we have
\bea
\frac{\langle \mu(z+\chi)+\gamma \rangle}{\langle \mu(z)+\gamma \rangle}=
\frac{1}{(q^{\gamma}t^{\mu}; q)_{\mu(\chi)}}, \quad 
\frac{\sigma(\mu(z+\chi))}{\sigma(\mu(z))}=(-1)^{\mu(\chi)}\prod_{\nu=0}^{\mu(\chi)-1}(q^{\nu}t^{\mu })^{-1}.
\ena
Hence, we can get a solution to (\ref{1}) in the following form:
\bea
\Phi(z)=t_{1}^{\alpha_{1}} \cdots t_{l}^{\alpha_{l}} \frac
        {\displaystyle \prod_{j=1}^{n'}\langle \mu_{j}'(z)+\gamma_{j}' \rangle}
        {\displaystyle \prod_{j=1}^{n}\langle \mu_{j}(z)+\gamma_{j} \rangle},
\label{5}
\ena
where $\alpha_{m}, \gamma_{j}, \gamma_{j}' \in {\Bbb C}\, , \mu_{j} , \mu_{j}' \in L^{*}$.

A function $\Phi(z)$ of type (\ref{5}) is called a $q$-{\it multiplicative function} at $|q|=1$.

\section{The $q$-twisted cohomology in Jordan-Pochhammer case}
Let us consider the $q$-multiplicative function of Jordan-Pochhammer type given by
\bea
\Phi(z)=t^{\alpha}\prod_{j=1}^{n}\frac{\langle z+\gamma_{j}' \rangle}{\langle z+\gamma_{j} \rangle},
\label{A1}
\ena
where $z \in {\Bbb C}\, , t=e^{2 \pi i \omega z}=q^{z}$ and $\alpha, \gamma_{j}, \gamma_{j}' \in {\Bbb C}\, $.
We assume $\gamma_{j} \not= \gamma_{k}'$ for any $j, k=1, \cdots , n$.

We denote by $D, D_{j}$ and $D_{j}'$ the difference operators corresponding to the displacements 
$\alpha \mapsto \alpha +1, \gamma_{j} \mapsto \gamma_{j}+1$ and $\gamma_{j}' \mapsto \gamma_{j}'+1$, respectively. 
Let $\cal A$ be the commutative algebra generated by $D^{\pm 1}, D_{j}^{\pm 1}$ 
and $D_{j}'^{\pm 1} \, (j=1, \cdots , n)$ over ${\Bbb C}\,$. 
We define a subspace $Z$ of ${\Bbb C}\, (t)$ by
\bea
Z:=\{(\kappa \Phi)/\Phi | \kappa \in {\cal A} \}.
\ena
It is easy to see that
\bea
Z=\left\{\frac{f(t)}{\prod_{j=1}^{n}(c_{j}q^{-\ell_{j}}t; q)_{\ell_{j}}(c_{j}'t; q)_{\ell_{j}'}}\Bigg|
    f(t) \in {\Bbb C}\, [t, t^{-1}] \, {\rm and} \, \ell_{j}, \ell_{j}' \in {\Bbb Z}_{\geq 0} \right\},
\ena
where $c_{j}=e^{2 \pi i \omega \gamma_{j}}$ and $c_{j}'=e^{2 \pi i \omega \gamma_{j}'}$. 

Next we define another space of rational functions. 
An important point is that the function $\langle x \rangle$ satisfies also the following functional relation:
\bea
\frac{\langle x+\frac{1}{\omega} \rangle}{\langle x \rangle}=\frac{1}{1-e^{2 \pi i x}}.
\label{5.1}
\ena
We denote by $\widetilde{D}, \widetilde{D_{j}}$ and $\widetilde{D_{j}'}$ the difference operators 
corresponding to the displacements 
$\alpha \mapsto \alpha+\frac{1}{\omega}, \gamma_{j} \mapsto \gamma_{j}+\frac{1}{\omega}$ and
$\gamma_{j}' \mapsto \gamma_{j}'+\frac{1}{\omega}$, respectively.
In the same way as before, we consider the commutative algebra $\widetilde{\cal A}$ 
generated by $\widetilde{D}^{\pm 1}, \widetilde{D_{j}}^{\pm 1}$ 
and $\widetilde{D_{j}'}^{\pm 1} \, (j=1, \cdots , n)$ over ${\Bbb C}\, $ 
and set
\bea
\widetilde{Z}:=\{(\widetilde{\kappa}\Phi)/\Phi | \widetilde{\kappa} \in \widetilde{\cal A} \}.
\ena
Then we have
\bea
\widetilde{Z}=
\left\{\frac{\widetilde{f}(T)}{\prod_{j=1}^{n}(C_{j}Q^{-\widetilde{\ell_{j}}}T; Q)_{\widetilde{\ell_{j}}}
                                              (C_{j}'T; Q)_{\widetilde{\ell_{j}'}}} 
 \Bigg|\widetilde{f}(T) \in {\Bbb C}\, [T, T^{-1}]\, 
  {\rm and } \, \widetilde{\ell_{j}}, \widetilde{\ell_{j}'} \in {\Bbb Z}_{\geq 0}\right\},
\ena
where $Q=e^{\frac{2 \pi i}{\omega}}, T=e^{2 \pi i z}, C_{j}=e^{2 \pi i \gamma_{j}}$ 
and $C_{j}'=e^{2 \pi i \gamma_{j}'}$.

Now we define a pairing between $Z$ and $\widetilde{Z}$. 
For
\bea
\varphi(t)=\frac{t^{m}}{\prod_{j=1}^{n}(c_{j}q^{-\ell_{j}}t; q)_{\ell_{j}}(c_{j}'t; q)_{\ell_{j}'}}
\in Z
\ena
and
\bea
\widetilde{\varphi}(T)=
\frac{T^{\widetilde{m}}}{\prod_{j=1}^{n}(C_{j}Q^{-\widetilde{\ell_{j}}}T; Q)_{\widetilde{\ell_{j}}}
                                       (C_{j}'T; Q)_{\widetilde{\ell_{j}'}}} 
\in \widetilde{Z},
\ena
we set
\bea
I(\varphi, \widetilde{\varphi}):=\int_{C}dz\Phi(z)\varphi(t)\widetilde{\varphi}(T).
\label{6}
\ena
Here the contour $C$ is taken to be the imaginary axis $(-i\infty, i\infty)$ except that the poles at
\bea
-\gamma_{j}+\ell_{j}+\frac{\widetilde{\ell_{j}}}{\omega}+{\Bbb Z}_{\leq 0}+\frac{1}{\omega}{\Bbb Z}_{\leq 0}
\quad (j=1, \cdots , n)
\ena
are on the left of $C$ and the poles at
\bea
-\gamma_{j}'-\ell_{j}'-\frac{\widetilde{\ell_{j}'}}{\omega}+{\Bbb Z}_{\geq 1}+\frac{1}{\omega}{\Bbb Z}_{\geq 1}
\quad (j=1, \cdots , n)
\ena
are on the right of $C$. 
Then the integral (\ref{6}) is absolutely convergent if
\bea
0< {\rm Re}\alpha +m+\frac{\widetilde{m}}{\omega}
<\sum_{j=1}^{n}({\rm Re}\gamma_{j}-{\rm Re}\gamma_{j}')+
 \sum_{j=1}^{n}(\ell_{j}+\ell_{j}')+\frac{1}{\omega}\sum_{j=1}^{n}(\widetilde{\ell_{j}}+\widetilde{\ell_{j}'}).
\label{7}
\ena
Under the condition (\ref{7}) the integrand of (\ref{6}) decreases exponentially as $z \to \pm i \infty$.

Let us consider cohomologies of $Z$ and $\widetilde{Z}$ associated with the integral (\ref{6}). 
We set
\bea
B:={\rm span}_{\Bbb C}\left\{\psi(t)-b_{\chi}(t)\psi(q^{\chi}t) | \psi(t) \in Z, \chi \in {\Bbb Z} \right\},
\label{7.1}
\ena
where $b_{\chi}(t)=\Phi(z+\chi)/\Phi(z)$. 
Note that for any $\psi(t) \in Z, \widetilde{\varphi}(T) \in \widetilde{Z}$ and $\chi \in {\Bbb Z}$ 
we can deform the contour $C$ so that there are no poles of the function
$\Phi(z)\psi(t)\widetilde{\varphi}(T)$ between $C$ and $C+\chi$.  
Thus we have
\bea
\int_{C}dz \Phi(z)\{\psi(t)-b_{\chi}(t)\psi(q^{\chi}t)\}\widetilde{\varphi}(T)=
\left( \int_{C}-\int_{C+\chi} \right) dz \Phi(z)\psi(t)\widetilde{\varphi}(T)=0,
\ena
if all the integrals are convergent.
Here we used the fact that $T=e^{2 \pi i z}$ is invariant under the change $z \to z-\chi$. 
Hence, we find
\bea
I(\varphi_{0}, \widetilde{\varphi})=0 \quad 
{\rm for} \, \varphi_{0} \in B \, {\rm and} \, \widetilde{\varphi} \in \widetilde{Z}.
\ena
Similarly, we set
\bea
\widetilde{B}:={\rm span}_{\Bbb C}\left\{\widetilde{\psi}(T)-\widetilde{b_{\chi}}(T)\widetilde{\psi}(Q^{\chi}T) |
                                         \widetilde{\psi}(T) \in \widetilde{Z}, \chi \in {\Bbb Z} \right\},
\ena
where $\widetilde{b_{\chi}}(T):=\Phi(z+\frac{\chi}{\omega})/\Phi(z)$. 
Then we have
\bea
I(\varphi, \widetilde{\varphi_{0}})=0 \quad 
{\rm for} \, \varphi \in Z \, {\rm and} \, \widetilde{\varphi_{0}} \in \widetilde{B}.
\ena
{}From these relations, we define the {\it $q$-twisted cohomology} $H$ and $\widetilde{H}$ by
\bea
H:=Z/B \quad {\rm and} \quad \widetilde{H}:=\widetilde{Z}/\widetilde{B}.
\ena

We note that the structure of the cohomology $H$ is determined by the parameters
$q, q^{\alpha}, c_{j}$ and $c_{j}' \, (j=1, \cdots , n)$.
We write down this dependence explicitly as
\bea
H={\rm H}(q| q^{\alpha}; c_{1}, \cdots , c_{n}; c_{1}', \cdots , c_{n}').
\ena
Then $\widetilde{H}$ is written as
\bea
\widetilde{H}={\rm H}(Q| Q^{\omega \alpha}; C_{1}, \cdots , C_{n}; C_{1}', \cdots , C_{n}').
\ena
It is easy to see that the following proposition holds.

\begin{prop}
The cohomology ${\rm H}(q| q^{\alpha}; c_{1}, \cdots , c_{n}; c_{1}', \cdots , c_{n}')$ is generated by
\bea
\left\{ \frac{1}{1-c_{j}'t} \Big| j=1, \cdots , n \right\},
\ena
if the parameters $q, q^{\alpha}, c_{j}$ and $c_{j}' \, (j=1, \cdots , n)$ are generic.
\end{prop}   

Moreover, we see that the set $\{\frac{1}{1-c_{j}'t}\}_{j=1, \cdots , n}$ is a basis of $H$ 
{}from the following determinant formula. 

\begin{prop}
Set
\bea
\varphi_{j}(t)=\frac{1}{1-c_{j}'t}, \quad \widetilde{\varphi}_{j}(T)=\frac{1}{1-C_{j}'T}, \quad 
(j=1, \cdots , n).
\ena
Then
\bea
\det{(I(\varphi_{j}, \widetilde{\varphi}_{k}))_{j, k=1, \cdots , n}}
&=& \langle 1 \rangle^{n} \exp{(-2 \pi i \omega \alpha \sum_{j=1}^{n}\gamma_{j}' )} \label{7.9}  \\
&\times& \frac{\langle \alpha+\sum_{j=1}^{n}\gamma_{j}-\sum_{j=1}^{n}\gamma_{j}' \rangle}
                  {\langle \alpha \rangle \prod_{j, k=1}^{n}\langle \gamma_{j}-\gamma_{k}' \rangle}\no \\ 
&\times& \prod_{1 \leq j<k \leq n}
             \frac{(1-e^{2 \pi i \omega (\gamma_{j}'-\gamma_{k}')})(1-e^{2 \pi i (\gamma_{j}'-\gamma_{k}')})}
                  {\sigma (\gamma_{j}'-\gamma_{k}')}. \no
\ena
\end{prop}

\proof
First we set
\bea
\quad
\psi_{k}(t)=\frac{1}{1-c_{k}'t}\prod_{j=1}^{k-1}\frac{1-c_{j}t}{1-c_{j}'t}, \quad
\widetilde{\psi}_{k}(T)=\frac{1}{1-C_{k}'T}\prod_{j=1}^{k-1}\frac{1-C_{j}T}{1-C_{j}'T}, \quad
(j=1, \cdots , n).
\ena
Then we have
\bea
\psi_{k}(t)=\sum_{p=1}^{k-1}\frac{1-c_{p}/c_{p}'}{1-c_{k}'/c_{p}'}
            \prod_{j=1 \atop j \not= p}^{k-1}\frac{1-c_{j}/c_{p}'}{1-c_{j}'/c_{p}'}\varphi_{p}(t)+
            \prod_{j=1}^{k-1}\frac{1-c_{j}/c_{k}'}{1-c_{j}'/c_{k}'}\varphi_{k}(t)
\ena
and a similar formula for $\widetilde{\psi}_{k}(T)$.

Hence we find
\bea
\det{(I(\varphi_{j}, \widetilde{\varphi}_{k}))}=
\prod_{1 \leq j<k \leq n}\left( \frac{1-c_{j}'/c_{k}'}{1-c_{j}/c_{k}'}\frac{1-C_{j}'/C_{k}'}{1-C_{j}/C_{k}'}
                         \right)
\det{(I(\psi_{j}, \widetilde{\psi}_{k}))}.
\label{8}
\ena
The determinant in the right hand side of (\ref{8}) is a special case of the determinant discussed in \cite{MT}.
Combining the result in \cite{MT} and (\ref{8}), we get the formula (\ref{7.9}).
\qed
\newline

\rem
In the case of $n=1$ and $\gamma_{1}'=0$, the formula (\ref{7.9}) is represented as follows:
\bea
\int_{C} q^{\alpha z}\frac{\langle z \rangle}{\langle z+\beta \rangle}\frac{1}{1-t}\frac{1}{1-T} dz=
\int_{C} q^{\alpha z}\frac{\langle z+1+\frac{1}{\omega} \rangle}{\langle z+\beta \rangle} dz=
\frac{\langle 1 \rangle \langle \alpha+\beta\rangle}{\langle \alpha \rangle \langle \beta \rangle}.
\label{9}
\ena
We may call (\ref{9}) the $q$-{\it Beta integral formula} at $|q|=1$.
\newline

To finish this section we find a system of difference equation in $\alpha$ satisfied by the function
\bea
\Psi(\alpha)=\Psi(\alpha | \widetilde{\varphi}):=\int_{C}dz \Phi(z)\widetilde{\varphi}(T)
\quad {\rm for} \, \widetilde{\varphi}(T) \in \widetilde{H}
\ena
in a similar manner to \cite{AK}.

For the $q$-multiplicative function (\ref{A1}), 
we can represent the funcion $b_{\chi}(t)=\frac{\Phi(z+\chi)}{\Phi(z)} \, (\chi \in {\Bbb Z})$ as follows:
\bea
b_{\chi}(t)=q^{\chi \alpha}\frac{b_{\chi}^{+}(t)}{b_{\chi}^{-}(t)},  
\ena
where $b_{\chi}^{+}(t)$ and $b_{\chi}^{-}(t)$ are polynomials in $t$ and have no common factor. 
For example, if $\chi =1$ we have
\bea
b_{1}^{+}(t)=\prod_{j=1}^{n}(1-c_{j}t), \quad b_{1}^{-}(t)=\prod_{j=1}^{n}(1-c_{j}'t).
\ena
By setting $\psi(t)=b_{\chi}^{-}(q^{-\chi}t)$ in (\ref{7.1}), we find
\bea
b_{\chi}^{-}(q^{-\chi}t)-q^{\chi \alpha}b_{\chi}^{+}(t) \in B.
\ena
Note that $t \, \Phi =D\Phi$, where $D$ is the difference operator defined by $\alpha \mapsto \alpha +1$. 
Therefore, we get
\bea
\left\{ b_{\chi}^{-}(q^{-\chi}D)-q^{\chi \alpha}b_{\chi}^{+}(D) \right\}\Psi =0 
\label{10}
\ena
for $\chi \in {\Bbb Z}$ such that $b_{\chi}^{-}(q^{-\chi}D)\Psi $ and $b_{\chi}^{+}(D)\Psi $ are defined.
These equations are $q$-analogues of Mellin-Sato hypergeometric equations \cite{AK} at $|q|=1$. 
In the case of $\chi =1$, the equation (\ref{10}) is given by
\bea
\left\{ \prod_{j=1}^{n}(1-q^{-1}c_{j}'D)-q^{\alpha}\prod_{j=1}^{n}(1-c_{j}D) \right\} \Psi =0.
\label{11}
\ena

\section{Application to the $q$-hypergeometric function}

\subsection{Preliminaries}
Following \cite{GR}, we recall some properties of 
the basic hypergeometric series with $0<|q|<1$ given by
\bea
\phi(a, b, c; t):=\sum_{k=0}^{\infty} 
\frac{(a; q)_{k}(b; q)_{k}}{(q; q)_{k}(c; q)_{k}}t^{k}
\ena
for $|t|<1$. 

This function satisfies the Heine's relations:
\bea
& & \phi(a, b, q^{-1}c)-\phi(a, b, c)=
       tc\frac{(1-a)(1-b)}{(q-c)(1-c)}\phi(qa, qb, qc),  \\
& & \phi(qa, b, c)-\phi(a, b, c)=
       ta\frac{1-b}{1-c}\phi(qa, qb, qc),  \\
& & \phi(qa, q^{-1}b, c)-\phi(a, b, c)=
       q^{-1}t\frac{aq-b}{1-c}\phi(qa, b, qc). 
\ena
Here we abbreviated $\phi(a, b, c; x)$ to $\phi(a, b, c)$. 

It also satisfies a connection formula:
\bea
\phi(a, b, c; t)&=&\frac{(b)_{\infty}(c/a)_{\infty}}{(c)_{\infty}(b/a)_{\infty}}
                   \frac{\Theta(at)}{\Theta(t)}\phi(a, qa/c, qa/b; qc/abt) \\
&+&\frac{(a)_{\infty}(c/b)_{\infty}}{(c)_{\infty}(a/b)_{\infty}}
   \frac{\Theta(bt)}{\Theta(t)}\phi(b, qb/c, qb/a; qc/abt), \no
\ena
where $(a)_{\infty}:=\prod_{j=1}^{\infty}(1-q^{j-1}a)$ and 
$\Theta(x):=(q)_{\infty}(x)_{\infty}(q/x)_{\infty}$. 
\newline

Now we consider the $q$-hypergeometric function of the Barnes type at $|q|=1$ \cite{NU}, 
which is defined as follows in our notation:
\bea
\Psi(\alpha, \beta, \gamma; x):=
\frac{\langle \alpha \rangle \langle \beta \rangle}{\langle 1 \rangle \langle \gamma \rangle}
\left( -\frac{1}{2 \pi i} \right)\int_{C_{0}}\frac{\langle z+1 \rangle \langle z+\gamma \rangle}
                                              {\langle z+\alpha \rangle \langle z+\beta \rangle}
                                         \frac{\pi (-q^{x})^{z}}{\sin{\pi z}} dz, 
\label{12}
\ena
where $-q^{x}=e^{2 \pi i \omega x+\pi i}$ and the contour $C_{0}$ is the imaginary axis $(-i \infty, i\infty)$ 
except that the poles at
\bea
-\alpha +{\Bbb Z}_{\leq 0}+\frac{1}{\omega}{\Bbb Z}_{\leq 0}, \quad
-\beta +{\Bbb Z}_{\leq 0}+\frac{1}{\omega}{\Bbb Z}_{\leq 0} 
\ena
are on the left of $C_{0}$ and the poles at
\bea
{\Bbb Z}_{\geq 0}+\frac{1}{\omega}{\Bbb Z}_{\geq 0}, \quad 
-\gamma +{\Bbb Z}_{\geq 1}+\frac{1}{\omega}{\Bbb Z}_{\geq 1}
\ena
are on the right of $C_{0}$.

By using 
\bea
\left( -\frac{1}{2 \pi i} \right) \frac{\pi (-q^{x})^{z}}{\sin{\pi z}}=q^{xz}\frac{1}{1-e^{2 \pi i z}},
\ena
we can rewrite (\ref{12}) as follows:
\bea
\Psi(\alpha, \beta, \gamma; x)=
\frac{\langle \alpha \rangle \langle \beta \rangle}{\langle 1 \rangle \langle \gamma \rangle}
\int_{C_{0}}q^{xz}\frac{\langle z+1+\frac{1}{\omega} \rangle \langle z+\gamma \rangle}
                   {\langle z+\alpha \rangle \langle z+\beta \rangle}dz.
\label{13}
\ena
Now we denote by $\Phi(z)$ the integrand of (\ref{13}):
\bea
\Phi(z)=q^{xz}\frac{\langle z+1+\frac{1}{\omega} \rangle \langle z+\gamma \rangle}
                   {\langle z+\alpha \rangle \langle z+\beta \rangle}.
\ena
This function $\Phi(z)$ is the $q$-multiplicative function of Jordan-Pochhammer type.
{}From (\ref{7}), the integral (\ref{13}) is absolutely convergent if
\bea
0< {\rm Re} x < 1+\frac{1}{\omega}+{\rm Re} \gamma -{\rm Re} \alpha -{\rm Re} \beta.
\label{14}
\ena
In this case, the equation (\ref{11}) is nothing but the hypergeometric difference equation at $|q|=1$:
\bea
\left\{ (1-D)(1-q^{\gamma -1}D)-q^{x}(1-q^{\alpha}D)(1-q^{\beta}D) \right\}\Psi =0,
\label{A1.5}
\ena
where $D$ is the difference operator defined by $x \mapsto x+1$.

For the $q$-multiplicative function (\ref{13}), 
we define two functional spaces $Z$ and $\widetilde{Z}$ as in the previous section, and set
\bea
\Psi(\alpha, \beta, \gamma; x|\widetilde{\varphi}):=
\frac{\langle \alpha \rangle \langle \beta \rangle}{\langle 1 \rangle \langle \gamma \rangle}
\int_{C_{0}}\Phi(z)\widetilde{\varphi}(T)dz
\qquad {\rm for} \quad \widetilde{\varphi} \in \widetilde{Z}.
\label{A2}
\ena
Note that $\Psi(\alpha, \beta, \gamma; x)=\Psi(\alpha, \beta, \gamma; x|1)$. 
Then the function (\ref{A2}) also satisfies (\ref{A1.5}).

For simplicity's sake, hereafter we use the following notation:
\bea
t=e^{2 \pi i \omega z}, \, a=e^{2 \pi i \omega \alpha}, \, b=e^{2 \pi i \omega \beta}, \, c=e^{2 \pi i \omega \gamma}
\label{14.1}
\ena
and
\bea
T=e^{2 \pi i z}, \, A=e^{2 \pi i \alpha}, \, B=e^{2 \pi i \beta}, \, C=e^{2 \pi i \gamma}.
\ena

\subsection{Heine's relations}
\begin{prop}
We abbreviate $\Psi(\alpha, \beta, \gamma; x)$ to $\Psi(\alpha, \beta, \gamma)$. 
Then the following equalities hold:
\bea
& & \Psi(\alpha, \beta, \gamma-1)-\Psi(\alpha, \beta, \gamma)=
       q^{x}c\frac{(1-a)(1-b)}{(q-c)(1-c)}\Psi(\alpha+1, \beta+1, \gamma+1), \label{14.8}\\
& & \Psi(\alpha+1, \beta, \gamma)-\Psi(\alpha, \beta, \gamma)=
       q^{x}a\frac{1-b}{1-c}\Psi(\alpha+1, \beta+1, \gamma+1), \label{14.85}\\
& & \Psi(\alpha+1, \beta-1, \gamma)-\Psi(\alpha, \beta, \gamma)=
       q^{x-1}\frac{aq-b}{1-c}\Psi(\alpha+1, \beta, \gamma+1). \label{14.9}
\ena
\end{prop}

\proof
For $f(t) \in Z$, we set
\bea
[f]:=\frac{\langle \alpha \rangle \langle \beta \rangle}{\langle 1 \rangle \langle \gamma \rangle}
     \int_{C_{0}}\Phi(z)f(t)dz.
\ena

First we prove (\ref{14.8}). It is easy to see that
\bea
\Psi(\alpha, \beta, \gamma-1)=\left[\frac{1-q^{-1}ct}{1-q^{-1}c}\right], \quad \Psi(\alpha, \beta, \gamma)=[1].
\ena
Hence, we have
\bea
\Psi(\alpha, \beta, \gamma-1)-\Psi(\alpha, \beta, \gamma)=\left[\frac{c(1-t)}{q-c}\right].
\label{15}
\ena
On the other hand, by changing the variable $z \to z-1$, we find
\bea
\Psi(\alpha+1, \beta+1, \gamma+1)=\left[q^{-x}\frac{(1-c)(1-t)}{(1-a)(1-b)}\right].
\label{16}
\ena
{}From (\ref{15}) and (\ref{16}), we get (\ref{14.8}).
We can prove (\ref{14.85}) in the same way as above.

Next we prove (\ref{14.9}). 
By changing the variable $z \to z+1$, we have
\bea
\Psi(\alpha+1, \beta-1, \gamma)=\left[ q^{x-1}\frac{(q-b)(1-at)(1-qat)}{(1-a)(1-qt)(1-ct)} \right].
\ena
It is easy to see that
\bea
\Psi(\alpha, \beta, \gamma)=[1], \quad 
\Psi(\alpha+1, \beta, \gamma+1)=\left[ \frac{(1-c)(1-at)}{(1-a)(1-ct)} \right].
\ena
By using this, we can find the following:
\bea
& & \Psi(\alpha+1, \beta-1, \gamma)-\Psi(\alpha, \beta, \gamma)-
    q^{x-1}\frac{aq-b}{1-c}\Psi(\alpha+1, \beta, \gamma+1) \label{17} \\
& & \qquad {}=\left[ -1+q^{x}\frac{(1-at)(1-bt)}{(1-qt)(1-ct)} \right]. \no
\ena
Note that
\bea
-1+q^{x}\frac{(1-at)(1-bt)}{(1-qt)(1-ct)}=-\{1-b_{1}(t)\cdot 1 \} \in B.
\ena
Therefore, (\ref{17}) equals to $0$. This completes the proof of (\ref{14.9}).
\qed\newline

In the proof above, we see that Heine's relations come from some relations in $H$. 
Hence, we find that the funcion $\Psi(\alpha, \beta, \gamma; x| \widetilde{\varphi})$ 
(\ref{A2}) also satisfies Heine's relations.

\subsection{Connection formula}
\begin{prop}
\bea
& & \Psi(\alpha, \beta, \gamma; x) \label{18} \\
& & {}=
\frac{\langle \beta \rangle \langle \gamma-\alpha \rangle}{\langle \gamma \rangle \langle \beta-\alpha \rangle}
\frac{\sigma(x+\alpha)}{\sigma(x)}
\Psi(\alpha, 1+\alpha-\gamma, 1+\alpha-\beta; 1+\frac{1}{\omega}+\gamma-\alpha-\beta-x) \no \\
& & {}+
\frac{\langle \alpha \rangle \langle \gamma-\beta \rangle}{\langle \gamma \rangle \langle \alpha-\beta \rangle}
\frac{\sigma(x+\beta)}{\sigma(x)}
\Psi(\beta, 1+\beta-\gamma, 1+\beta-\alpha; 1+\frac{1}{\omega}+\gamma-\alpha-\beta-x).
\no
\ena
\end{prop}
\proof
We rewrite the integral
\bea
& & \Psi(\alpha, 1+\alpha-\gamma, 1+\alpha-\beta; 1+\frac{1}{\omega}+\gamma-\alpha-\beta-x) \label{19} \\
& & \qquad {}=\frac{\langle \alpha \rangle}{\langle 1 \rangle} \int_{C_{0}'} 
       q^{(1+\frac{1}{\omega}+\gamma-\alpha-\beta-x)z}
       \frac{\langle 1+\alpha-\gamma \rangle}{\langle 1+\alpha-\beta \rangle}
       \frac{\langle z+1+\frac{1}{\omega} \rangle \langle z+1+\alpha-\beta \rangle}
            {\langle z+\alpha \rangle \langle z+1+\alpha-\gamma \rangle}dz, \no
\ena
where $C_{0}'$ is the contour associated with the set of parameters $(\alpha, 1+\alpha-\gamma, 1+\alpha-\beta)$. 

By changing the variable $z \to -z-\alpha$, we have
\bea
{}\qquad
(\ref{19})=\frac{\langle \alpha \rangle}{\langle 1 \rangle}\!\int_{C_{0}}\! 
       q^{(1+\frac{1}{\omega}+\gamma-\alpha-\beta-x)(-z-\alpha)}
       \frac{\langle 1+\alpha-\gamma \rangle}{\langle 1+\alpha-\beta \rangle}
       \frac{\langle -z-\alpha+1+\frac{1}{\omega} \rangle \langle -z+1-\beta \rangle}
            {\langle -z \rangle \langle -z+1-\gamma \rangle}dz,
\label{20}
\ena
where $C_{0}$ is the contour defined in (\ref{12}). 
By using (\ref{4.1}), we have
\bea
& & {\rm the \, integrand \, of \, (\ref{20})}  \\
& & \qquad {}=q^{(1+\frac{1}{\omega}+\gamma-\alpha-\beta-x)(-z-\alpha)}
       \frac{\sigma(1+\alpha-\gamma)}{\sigma(1+\alpha-\beta)}
       \frac{\sigma(-z-\alpha+1+\frac{1}{\omega})\sigma(-z+1-\beta)}
            {\sigma(-z)\sigma(-z+1-\gamma)}\no \\
& & \qquad {}\times \frac{\langle \frac{1}{\omega}+\beta-\alpha \rangle}{\langle \frac{1}{\omega}+\gamma-\alpha \rangle}
             \frac{\langle z+1+\frac{1}{\omega} \rangle \langle z+\gamma+\frac{1}{\omega} \rangle}
                  {\langle z+\alpha \rangle \langle z+\beta+\frac{1}{\omega} \rangle}. \no
\ena
It can be shown that
\bea
& &
q^{(1+\frac{1}{\omega}+\gamma-\alpha-\beta-x)(-z-\alpha)}
\frac{\sigma(1+\alpha-\gamma)}{\sigma(1+\alpha-\beta)}
\frac{\sigma(-z-\alpha+1+\frac{1}{\omega})\sigma(-z+1-\beta)}
     {\sigma(-z)\sigma(-z+1-\gamma)} \label{21}\\
& & \quad {}=q^{xz}\frac{\sigma(x)}{\sigma(x+\alpha)}. \no
\ena
{}From (\ref{5.1}), we have
\bea
& & 
\frac{\langle \frac{1}{\omega}+\beta-\alpha \rangle}{\langle \frac{1}{\omega}+\gamma-\alpha \rangle}
\frac{\langle z+1+\frac{1}{\omega} \rangle \langle z+\gamma+\frac{1}{\omega} \rangle}
     {\langle z+\alpha \rangle \langle z+\beta+\frac{1}{\omega} \rangle} \label{22} \\
& &\quad {}=
\frac{\langle \beta-\alpha \rangle}{\langle \gamma-\alpha \rangle}
\frac{\langle z+1+\frac{1}{\omega} \rangle \langle z+\gamma \rangle}
     {\langle z+\alpha \rangle \langle z+\beta \rangle}
\frac{(A-C)(1-BT)}{(A-B)(1-CT)}. \no
\ena
Combining (\ref{21}) and (\ref{22}), we get
\bea
& & \Psi(\alpha, 1+\alpha-\gamma, 1+\alpha-\beta; 1+\frac{1}{\omega}+\gamma-\alpha-\beta-x)  \\
& & \qquad {}=
\frac{\langle \alpha \rangle \langle \beta-\alpha \rangle}{\langle 1 \rangle \langle \gamma-\alpha \rangle}
\frac{\sigma(x)}{\sigma(x+\alpha)} \int_{C_{0}}\Phi(z)\frac{(A-C)(1-BT)}{(A-B)(1-CT)}dz. \no
\ena

By exchanging $\alpha$ and $\beta$, we find
\bea
& & \Psi(\beta, 1+\beta-\gamma, 1+\beta-\alpha; 1+\frac{1}{\omega}+\gamma-\alpha-\beta-x)   \\
& & \qquad {}=
\frac{\langle \beta \rangle \langle \alpha-\beta \rangle}{\langle 1 \rangle \langle \gamma-\beta \rangle}
\frac{\sigma(x)}{\sigma(x+\beta)} \int_{C_{0}}\Phi(z)\frac{(B-C)(1-AT)}{(B-A)(1-CT)}dz. \no
\ena
Therefore, we get
\bea
{}\quad 
{\rm the \, rhs \, of \, (\ref{18})}&=&
\frac{\langle \alpha \rangle \langle \beta \rangle}{\langle 1 \rangle \langle \gamma \rangle}
\int_{C_{0}}\Phi(z)\left\{ \frac{(A-C)(1-BT)}{(A-B)(1-CT)}+\frac{(B-C)(1-AT)}{(B-A)(1-CT)} \right\}dz  \\ 
&=& 
\frac{\langle \alpha \rangle \langle \beta \rangle}{\langle 1 \rangle \langle \gamma \rangle}
\int_{C_{0}}\Phi(z) \cdot 1 dz=\Psi(\alpha, \beta, \gamma; x). \no
\ena
\qed\newline
 
In the proof above, we see that the formula (\ref{18}) comes from the following simple relation in $\widetilde{H}$:
\bea
\frac{(A-C)(1-BT)}{(A-B)(1-CT)}+\frac{(B-C)(1-AT)}{(B-A)(1-CT)}=1.
\label{23}
\ena
If we consider $H$ as a cohomology and $\widetilde{H}$ as its dual, that is a homology, 
then the relation (\ref{23}) is a relation among some homologies, 
and the formula (\ref{18}) is a linear relation among the integrals 
associated with different homologies.

\section{The $q$-little Jacobi polynomials at $|q|=1$}
First we recall the definition of the $q$-little Jacobi polynomials 
in the case of $0<|q|<1$ \cite{GR}:
\bea
p_{n}^{(\alpha, \beta)}(t):=\phi(q^{-n}, q^{\alpha+\beta+n+1}, q^{\alpha+1}; qt), \quad (n=0, 1, \cdots ).
\ena
The following orthogonality relation holds \cite{AA, GR}:
\bea
\int_{0}^{1}t^{\alpha-1}\frac{(tq)_{\infty}}{(tq^{\beta})_{\infty}}
p_{m}^{(\alpha-1, \beta-1)}(t)p_{n}^{(\alpha-1, \beta-1)}(t)d_{q}t=
\delta_{m, n}c_{n}, 
\label{24}
\ena
where
\bea
c_{n}=(1-q)\frac{(q)_{\infty}(q^{\alpha+\beta})_{\infty}}{(q^{\alpha})_{\infty}(q^{\beta})_{\infty}}
      \frac{1-q^{\alpha+\beta-1}}{1-q^{\alpha+\beta+2n-1}}
      \frac{(q)_{n}(q^{\beta})_{n}}{(q^{\alpha+\beta-1})_{n}(q^{\alpha})_{n}}q^{n\alpha}.
\ena
In (\ref{24}), the integral is a Jackson integral defined by
\bea
\int_{0}^{1}f(t)d_{q}t:=(1-q)\sum_{k=0}^{\infty}f(q^{n})q^{n}, 
\ena
and $(a)_{n}=(a; q)_{n}$.

The formula (\ref{24}) means that the $q$-little Jacobi polynomials are orthogonal polynomials 
with respect to the kernel of the $q$-Beta integral given by
\bea
\int_{0}^{1}t^{\alpha-1}\frac{(tq)_{\infty}}{(tq^{\beta})_{\infty}}d_{q}t=
(1-q)\frac{(q)_{\infty}(q^{\alpha+\beta})_{\infty}}{(q^{\alpha})_{\infty}(q^{\beta})_{\infty}}.
\label{24.1}
\ena

Let us consider the case of $|q|=1$. 
We can get the $q$-little Jacobi polynomials at $|q|=1$ from the $q$-hypergeometric function (\ref{13}) as follows.

\begin{prop}
For $n \in {\Bbb Z}_{\geq 0}$, we have
\bea
\lim_{\alpha \to -n}\Psi(\alpha, \beta, \gamma; x)=\phi(q^{-n}, q^{\beta}, q^{\gamma}; q^{x}).
\label{25}
\ena
Note that the right hand side of (\ref{25}) is a polynomial in $q^{x}$
and so makes sense at $|q|=1$.
\end{prop}

\proof
Recall the definition of $\Psi(\alpha, \beta, \gamma; x)$:
\bea
\Psi(\alpha, \beta, \gamma; x):=
\frac{\langle \alpha \rangle \langle \beta \rangle}{\langle 1 \rangle \langle \gamma \rangle}
\int_{C_{0}}\Phi(z)dz.
\ena
At $\alpha=-n$, the coefficient $\langle \alpha \rangle$ has a zero and
the integral has a pole because of pinches of the contour $C_{0}$ by poles at
$z=0, 1, \cdots , n$ and $z=-\alpha-n, -\alpha-n+1, \cdots , -\alpha$, respectively.
In order to avoid these pinches, we take the residues at $z=-\alpha-n, \cdots , -\alpha$. 
Then we get
\bea
\Psi(\alpha, \beta, \gamma; x)&=&
\frac{\langle \alpha \rangle \langle \beta \rangle}{\langle 1 \rangle \langle \gamma \rangle}
\sum_{k=0}^{n}2 \pi i {\rm res}_{z=-\alpha-k}\Phi(z)dz \label{26} \\
&+&  
\frac{\langle \alpha \rangle \langle \beta \rangle}{\langle 1 \rangle \langle \gamma \rangle}
\times({\rm regular \, at} \, \alpha=-n). \no
\ena
The second term of the rhs of (\ref{26}) equals zero at $\alpha =-n$. 
Hence it suffices to calculate the limit of the first term.

By using (\ref{4.9}), we have
\bea
& & 2 \pi i{\rm res}_{z=-\alpha-(n-k)}\Phi(z)dz  \\
& & {}=
2 \pi i{\rm res}_{z=-\alpha}\Phi(z)\prod_{j=1}^{n-k}
\frac{(1-q^{-j+1}t)(1-q^{-j}ct)}{(1-q^{-j}at)(1-q^{-j}bt)}q^{-(n-k)x}dz \no \\
& & {}=
\langle 1 \rangle \langle 1+\frac{1}{\omega}-\alpha \rangle 
\frac{\langle \gamma-\alpha \rangle}{\langle \beta-\alpha \rangle}
\prod_{j=1}^{n-k}\frac{(1-q^{-j+1}/a)(1-q^{-j}c/a)}{(1-q^{-j})(1-q^{-j}b/a)}q^{-(\alpha+n-k)x}. \no
\ena
Here we used the notation (\ref{14.1}) and 
\bea
2 \pi i {\rm res}_{z=0}\frac{dz}{\langle z \rangle}=\frac{i}{\sqrt{\omega}}=\langle 1 \rangle.
\ena

Therefore, we get
\bea
& & \lim_{\alpha \to -n}\Psi(\alpha, \beta, \gamma; x)  \\
& & {}=
\lim_{\alpha \to -n}\langle \alpha \rangle \langle 1+\frac{1}{\omega}-\alpha \rangle
\frac{\langle \beta \rangle}{\langle \beta-\alpha \rangle}
\frac{\langle \gamma-\alpha \rangle}{\langle \gamma \rangle}
\sum_{k=0}^{n}q^{-(\alpha+n-k)x}
\prod_{j=1}^{n-k}\frac{(1-q^{-j+1}/a)(1-q^{-j}c/a)}{(1-q^{-j})(1-q^{-j}b/a)} \no \\
& & {}=
\sigma(-n)
\frac{\langle \beta \rangle}{\langle \beta+n \rangle}
\frac{\langle \gamma+n \rangle}{\langle \gamma \rangle}
\sum_{k=0}^{n}q^{kx}\prod_{j=1}^{n-k}
\frac{(1-q^{n+1-j})(1-q^{n-j}c)}{(1-q^{-j})(1-q^{n-j}b)} \no \\
& & {}=
(-1)^{n}q^{-\frac{n(n+1)}{2}}\prod_{j=1}^{n}\frac{1-q^{j}}{1-q^{-j}}
\sum_{k=0}^{n}q^{kx}\prod_{j=0}^{k-1}
\frac{(1-q^{-n+j})(1-q^{j}b)}{(1-q^{j})(1-q^{j}c)} \no \\
& & {}=\phi(q^{-n}, q^{\beta}, q^{\gamma}; q^{x}). \no
\ena
\qed
\newline

{}From this proposition, we get
\bea
p_{n}^{(\alpha, \beta)}(q^{x})=\Psi(-n, \alpha+\beta+n+1, \alpha+1; x+1), \quad (n=0, 1, \cdots ).
\label{27}
\ena

Then we find that the $q$-little Jacobi polynomials (\ref{27}) satisfy 
the orthogonal relation associated with the $q$-Beta integral at $|q|=1$ (\ref{9}).

\begin{prop}
\bea
\int_{C}q^{\alpha z}\frac{\langle z+1+\frac{1}{\omega} \rangle}{\langle z+\beta \rangle}
p_{m}^{(\alpha -1, \beta -1)}(t)p_{n}^{(\alpha-1, \beta-1)}(t)dz
=\delta_{m, n}c_{n}, 
\label{28}
\ena
where $t=q^{z}=e^{2 \pi i \omega z}$ and
\bea
c_{n}=\frac{\langle 1 \rangle \langle \alpha+\beta \rangle}{\langle \alpha \rangle \langle \beta \rangle}
      \frac{1-q^{\alpha+\beta-1}}{1-q^{\alpha+\beta+2n-1}}
      \frac{(q)_{n}(q^{\beta})_{n}}{(q^{\alpha+\beta-1})_{n}(q^{\alpha})_{n}}q^{n\alpha}.
\ena
In the left hand side of (\ref{28}), the contour $C$ is the imaginary axis $(-i\infty, i\infty)$ 
except that the poles at ${\Bbb Z}_{\geq 0}+\frac{1}{\omega}{\Bbb Z}_{\geq 0}$ are on the right of $C$ 
and the poles at $-\beta+{\Bbb Z}_{\leq 0}+\frac{1}{\omega}{\Bbb Z}_{\leq 0}$ are on the left of $C$.
\end{prop}

\proof
First we rewrite the orthogonality relation with $0<|q|<1$ (\ref{24}) as follows.
We expand the product
\bea
p_{m}^{(\alpha-1, \beta-1)}(t)p_{n}^{(\alpha-1, \beta-1)}(t)=\sum_{k=0}^{m+n}A_{k}^{m, n}t^{k}, 
\quad A_{k}^{m, n} \in {\Bbb C} \, .
\ena
By using (\ref{24.1}), we have
\bea
& & \int_{0}^{1}t^{\alpha-1}\frac{(tq)_{\infty}}{(tq^{\beta})_{\infty}}
p_{m}^{(\alpha-1, \beta-1)}(t)p_{n}^{(\alpha-1, \beta-1)}(t)d_{q}t \label{28.9} \\
& & \quad {}=\sum_{k=0}^{m+n} A_{k}^{m, n} \int_{0}^{1}t^{\alpha+k-1}\frac{(tq)_{\infty}}{(tq^{\beta})_{\infty}}d_{q}t \no \\
& & \quad {}=\sum_{k=0}^{m+n}A_{k}^{m, n}
    (1-q)\frac{(q)_{\infty}(q^{\alpha+\beta+k})_{\infty}}{(q^{\alpha+k-1})_{\infty}(q^{\beta})_{\infty}} \no \\
& & \quad {}=(1-q)\frac{(q)_{\infty}(q^{\alpha+\beta})_{\infty}}{(q^{\alpha})_{\infty}(q^{\beta})_{\infty}}
    \sum_{k=0}^{m+n}A_{k}^{m, n}\prod_{j=0}^{k-1}\frac{1-q^{\alpha+j}}{1-q^{\alpha+\beta+j}}. \no
\ena
Hence the relation (\ref{24}) is equivalent to
\bea
\sum_{k=0}^{m+n}A_{k}^{m, n}\prod_{j=0}^{k-1}\frac{1-q^{\alpha+j}}{1-q^{\alpha+\beta+j}}=
\delta_{m, n}\frac{1-q^{\alpha+\beta-1}}{1-q^{\alpha+\beta+2n-1}}
\frac{(q)_{n}(q^{\beta})_{n}}{(q^{\alpha+\beta-1})_{n}(q^{\alpha})_{n}}q^{n\alpha}.
\label{29}
\ena
Note that (\ref{29}) is an algebraic equality and holds also in the case of $|q|=1$. 

On the other hand, we find the following from (\ref{9}) in the same way as (\ref{28.9}):
\bea
{}\qquad 
\int_{C}q^{\alpha z}\frac{\langle z+1+\frac{1}{\omega} \rangle}{\langle z+\beta \rangle}
p_{m}^{(\alpha -1, \beta -1)}(t)p_{n}^{(\alpha-1, \beta-1)}(t)dz
=
\frac{\langle 1 \rangle \langle \alpha+\beta \rangle}{\langle \alpha \rangle \langle \beta \rangle}
\sum_{k=0}^{m+n}A_{k}^{m, n}\prod_{j=0}^{k-1}\frac{1-q^{\alpha+j}}{1-q^{\alpha+\beta+j}}.
\label{30}
\ena
{}From (\ref{29}) and (\ref{30}), we get (\ref{28}).
\qed

\section*{Acknoledgement}
The author thanks Tetsuji Miwa for valuable remarks. 

%


\begin{thebibliography}{II}
\parsep .05\bls
\itemsep 0pt
\frenchspacing

\bibitem[A1]{A1}
Aomoto, K\&, \itl{$q$-analogue of de Rham cohomology associated with Jackson integrals. 
{\rm I, II}}, Proc. Japan Acad., {\bf 66} Ser. A (1990), 161-164, 240-244.

\bibitem[A2]{A2}
Aomoto, K\&, \itl{Finiteness of a cohomology associated with certain Jackson integrals}, 
Tohoku Math. J., {\bf 43} (1991), 75-101.

\bibitem[AA]{AA}
Andrews, G. E. and Askey, R\&, \itl{Enumeration of partitions: the role of Eulerian series 
and $q$-orthogonal polynomials}, 
Higher Combinatorics, M. Aigner, eds., Reidel, Boston, Mass. (1977), 3-26. 

\bibitem[AK]{AK}
Aomoto, K. and Kato, Y\&, \itl{A $q$-analogue of de Rham cohomology associated with Jackson integrals}, 
Special Functions, M. Kashiwara and T. Miwa, eds., Proc. of the Hayashibara Forum, Springer, 1990.

\bibitem[B]{B}
Barnes, E. W\&, \itl{Theory of the double gamma functions}, 
Phil. Trans. Roy. Soc. A, {\bf 196} (1901), 265-388.

\bibitem[F]{F}
Faddeev, L\&, \itl{Modular double of quantum group}, math.QA/9912078, and references therein.

\bibitem[GR]{GR}
Gasper, G. and Rahman, M\&, \itl{Basic hypergeometric series}, 
Encyclopedia of Mathmatics and its Applications 35. Cambridge Univ. Press., 1990.

\bibitem[JM]{JM}
Jimbo, M. and Miwa, T\&, \itl{Quantum KZ equation with $|q|=1$ and correlation
functions of the XXZ model in the gapless regime}, J. Phys. A: Math. Gen.
{\bf 29} (1996), 2923-2958.

\bibitem[MT]{MT}
Miwa, T. and Takeyama, Y\&, \itl{Determinant formula for the solutions of
the quantum Knizhnik-Zamolodchikov equation with $|q|=1$}, Contemporary Mathmatics, {\bf 248} 
(1999), 377-393.

\bibitem[N]{N}
Nishizawa, M\&, \itl{On a solution of a $q$-difference analogue of 
Lauricella's $D$-type hypergeometric equation with $|q|=1$}, 
Publ. RIMS. Kyoto Univ., {\bf 34} (1998), 277-290.

\bibitem[NU]{NU}
Nishizawa, M. and Ueno, K\&, \itl{Integral solutions of $q$-difference equations 
of the hypergeometric type with $|q|=1$}, Proceedings of the workshop 
``Infinite Analysis - Integral Systems and Representation Theory'', IIAS Report No. 1997-001, 247-255.

\bibitem[Sa]{SSM}
Sato, M., Shintani, T. and Muro, M\&, \itl{Theory of prehomogeneous vector spaces (Algebraic part)}, 
Nagoya Math. J., {\bf 120} (1990), 1-33. 

\bibitem[Sm]{S}
Smirnov, F. A\&, \itl{Dual Baxter equations and quantization of affine Jacobian}, 
math-ph/0001032.

\bibitem[TV1]{TV1}
Tarasov, V. and Varchenko, A\&, \itl{Geometry of $q$-hypergeometric functions as
a bridge between Yangians and quantum affine algebras}, Invent. Math. {\bf 128}
(1997), 501--588. 

\bibitem[TV2]{TV2}
Tarasov, V. and Varchenko, A\&, \itl{Geometry of $q$-hypergeometric functions,
quantum affine algebras and elliptic quantum groups}, Ast\'erisque {\bf 246}
(1997), 1--135.

\end{thebibliography}
\end{document}